\documentclass{amsart}

\newtheorem{theorem}{Theorem}
\newtheorem{lemma}{Lemma}

\newcommand{\nwc}{\newcommand}
\nwc{\qref}[1]{(\ref{#1})} 

\renewcommand{\div}{\nabla\!\cdot\!}
\nwc{\udotgrad}{\uu\!\cdot\!\nabla}

\nwc{\dsl}{\displaystyle}

\newcommand{\PP}{{\mathcal P}}
\newcommand{\QQ}{{\mathcal Q}}

\nwc{\pcon}{\beta}  
\nwc{\pavs}{\bar{p}_s}
\nwc{\epss}{{\varepsilon}}
\nwc{\pgh}{p_{gh}}
\nwc{\SP}{{\cal S}_p}
\nwc{\SPb}{{\cal S}_\Gamma}
\nwc{\cut}{\xi}

\nwc{\Gdag}{G^{\dag}}
\nwc{\I}{{\cal I}}
\nwc{\gam}{\gamma}
\nwc{\yy}{\vec{y}}
\nwc{\zz}{\vec{z}}
\nwc{\hu}{\hat{u}}
\nwc{\hv}{\hat{v}}

\nwc{\xuh}{X_h}
\nwc{\xph}{Y_h}
\nwc{\Tmax}{T_{\rm max}}
\nwc{\Phifac}{\Phi\grad\nn}
\nwc{\tf}{\tilde{f}}
\nwc{\lapg}{\Delta_\G}
\nwc{\gradg}{\nabla_\G}
\nwc{\lapG}{\Delta_\Gamma}
\nwc{\lapI}{\Delta_\I}
\nwc{\gradG}{\nabla_\Gamma}
\nwc{\gradI}{\nabla_\I}
\nwc{\N}{\mathbb{N}}
\nwc{\ip}[1]{\langle #1 \rangle}
\nwc{\Omgs}{{\Omega_s}}
\nwc{\G}{{\cal G}}

\nwc{\ffd}{\ff_{\Dt}}
\nwc{\UU}{\vec{U}}
\nwc{\emb}{\hookrightarrow}
\nwc{\shalf}{{\textstyle\frac12}}

\nwc{\uspace}{{H^2\cap H^1_0(\Omega,\R^N)}}

\nwc{\uuin}{\uu_{\rm in}}
\nwc{\hin}{h_{\rm in}}

\nwc{\Htwoloc}{H^2_{\rm loc}}

\newcommand{\Hdiv}{H({\rm div};\Omega)}

\newcommand{\dist}{\mathop{\rm dist}\nolimits}
\newcommand{\lap}{\Delta }
\newcommand{\Omg}{\Omega}
\newcommand{\bOmg}{\bar{\Omg}}
\newcommand{\Gam}{\Gamma}
\newcommand{\Omgin}{\Omg \backslash \Omg_{s}}

\newcommand{\Omgbdry}{\Omg_{s}}

\newcommand{\eps}{\varepsilon}

\nwc{\ba}{\vec{a}}
\newcommand{\nn}{\vec{n}}

\newcommand{\ff}{\vec{f}}

\renewcommand{\gg}{\vec{g}}
\newcommand{\uu}{\vec{u}}
\newcommand{\vv}{\vec{v}}

\newcommand{\xx}{\vec{x}}

\newcommand{\nproj}{\nn\nn^t}
\newcommand{\tanproj}{I-\nproj}
\newcommand{\upe}{\vec{u}_{\perp}}
\newcommand{\upa}{\vec{u}_{\parallel}}

\newcommand{\grad}{\nabla}
\newcommand{\nder}{\nn\cdot\!\grad}
\newcommand{\pa}{\partial}

\newcommand{\R}{\mathbb{R}}

\renewcommand{\P}{\PP}

\newcommand{\Dt}{\Delta t}
\newcommand{\noind}{\noindent}
\newcommand{\pe}{p_{\mbox{\tiny E}}}
\newcommand{\ps}{p_{\mbox{\tiny S}}}

\renewcommand{\div}{\nabla \cdot}

\renewcommand{\(}{\big(}
\renewcommand{\)}{\big)}
\newcommand{\<}{\big<}
\renewcommand{\>}{\big>}

\usepackage{amsmath}
\usepackage{amsthm}
\usepackage{amssymb}

\usepackage{url}

\begin{document}

\title[On Incompressible Navier-Stokes Dynamics]{On Incompressible Navier-Stokes Dynamics: A New Approach for Analysis and Computation}

\author{Jian-Guo~Liu}
\address{Department of Mathematics \&
Institute for Physical Science and Technology,
University of Maryland \\ College Park MD 20742, USA.}
\email{jliu@math.umd.edu}

\author{Jie~Liu}
\address{Department of Mathematics,
University of Maryland \\ College Park MD 20742, USA.}
\email{jieliu@math.umd.edu}

\author{Robert~L.~Pego}
\address{Department of Mathematical Sciences, Carnegie Mellon University \\
Pittsburgh, PA 15213, USA.}
\email{rpego@cmu.edu}

\begin{abstract}
We show that in bounded domains with no-slip boundary
conditions, the Navier-Stokes pressure can be determined in a such
way that it is strictly dominated by viscosity. As a consequence,
in a general domain we can treat the Navier-Stokes equations as a
perturbed vector diffusion equation, instead of as a perturbed
Stokes system.
To illustrate the advantages of this view, we provide
a simple proof of the unconditional stability
of a difference scheme that is implicit only in viscosity and
explicit in both pressure and convection terms, requiring no
solution of stationary Stokes systems or inf-sup conditions. 
\end{abstract}

\maketitle

\section{Introduction} \label{S.introduction}
Consider the Navier-Stokes equations for incompressible fluid flow
in $\Omg$ with no-slip boundary conditions on $\Gamma:= \pa\Omg$:
\begin{align}
  \pa_t \uu + \udotgrad \uu + \nabla p &= \nu \Delta \uu + \ff
  \quad\mbox{  in $\Omg$},  \label{NSE1} \\
  \div \uu &= 0
  \quad\mbox{  in $\Omg$}, \label{NSE2} \\
  \quad \uu &=0  \quad\mbox{ on $\Gamma$}.  \label{NSE3}
\end{align}
Here $\uu$ is the fluid velocity, $p$ the pressure, and $\nu$ is
the kinematic viscosity coefficient, assumed to be a fixed
positive constant. Pressure plays a role like a Lagrange
multiplier to enforce the incompressibility constraint, and this
has been a main source of difficulties. In this paper we will describe
some of the main results of \cite{LLP}, which indicate
that the pressure can be obtained in a way that
leads to considerable simplifications in both computation and
analysis.

A standard way to determine $p$ is via the Helmholtz-Hodge
decomposition. We can rewrite \qref{NSE1} as
\begin{equation}
 \pa_t \uu + \PP( \udotgrad \uu - \ff - \nu \Delta \uu) = 0. \label{NSE4}
\end{equation}
where $\PP$ is the standard Helmholtz projection operator onto
divergence-free fields. 

In this formulation, solutions formally satisfy
$\pa_t(\div\uu)=0$. The dissipation in \qref{NSE4} appears
degenerate due to the fact that $\PP$ annihilates gradients, so
{the analysis of \qref{NSE4} is usually restricted to spaces
of divergence-free fields}. But alternatives are possible if the 
pressure is determined differently when the velocity field has
non-zero divergence. Instead of \qref{NSE4}, we will
consider the unconstrained formulation
\begin{equation}
 \pa_t \uu + \PP( \udotgrad \uu - \ff - \nu \Delta \uu )
 = \nu \nabla(\div \uu). \label{NSE5}
\end{equation}
There is no difference as long as $\div\uu=0$, of course. But we
argue that \qref{NSE5} enjoys superior stability properties, for
two reasons.  The first is heuristic. The incompressibility
constraint is enforced in a more robust way, because the
divergence of velocity satisfies a diffusion equation with no-flux
boundary conditions. Naturally, if $\div\uu = 0$ initially, this
remains true for all later time, and one has a solution of the
standard Navier-Stokes equations \qref{NSE1}--\qref{NSE3}.

The second reason is more profound. To explain, we recast
\qref{NSE5} in the form \qref{NSE1} while explicitly identifying
the separate contributions to the pressure term made by the
convection and viscosity terms. Using the Helmholtz projection
operator $\PP$, we introduce the {\em Euler pressure} $\pe$ and
{\em Stokes pressure} $\ps$ via the relations
\begin{align}
& \PP( \udotgrad \uu - \ff) = \udotgrad \uu - \ff +
\nabla \pe , \label{EulerP}
\\[4pt] 
& \PP(-\Delta \uu ) = -\Delta \uu
+\nabla(\div\uu) + \nabla \ps . \label{StokesP}
\end{align}
This puts \qref{NSE5} into the form \qref{NSE1} with $p=\pe+\nu\ps$:
\begin{equation}
 \pa_t \uu + \udotgrad \uu + \nabla \pe + \nu \nabla \ps
 = \nu \Delta \uu + \ff.  \label{NSE6}
\end{equation}
Since the Helmholtz projection is $L^2$-orthogonal, naturally the Stokes
pressure from \qref{StokesP} satisfies
\begin{equation}
\int_\Omg |\grad\ps|^2 \le \int_\Omg |\lap\uu|^2
\quad\mbox{if $\div\uu=0$.}
\end{equation}
In fact, it turns out that the Stokes pressure term
is {\em strictly} dominated by the viscosity term,
regardless of the divergence constraint.
The estimate contained in the following theorem from \cite{LLP}
is key to a new treatment of Navier-Stokes dynamics in bounded domains.
\begin{theorem} \label{T.main}
Let $\Omg \subset \R^N$ ($N\ge2$) be a connected bounded domain with $C^3$
boundary. Then for any $\eps>0$, there exists $C\ge0$ 
such that for all vector fields $\uu \in \uspace$,
the Stokes pressure $\ps$ determined by
\qref{StokesP} satisfies
\begin{equation}
 \int_\Omg |\nabla \ps |^2 \le 
 \pcon \int_\Omg |\lap
   \uu|^2 + C \int_\Omg |\nabla \uu|^2.
 \label{StokesPE}
\end{equation}
where $\pcon=\frac23+\eps$.
\end{theorem}

In this paper we will prove a slightly weaker version of this result
that has a simpler proof. Namely, we will show that the estimate
\qref{StokesPE} holds for {\em some} $\beta<1$.
The full result that \qref{StokesPE} holds with 
any $\beta>2/3$, independent of the domain or the space dimension,
is proved in \cite{LLP} using differential geometry to establish
sharp integrated Neumann-to-Dirichlet estimates for functions harmonic near
the boundary of $\Omg$.

Due to Theorem~\ref{T.main}, we can treat the Navier-Stokes
equations in bounded domains simply as a perturbation of the
vector diffusion equation $\pa_t\uu = \nu\Delta\uu$, regarding
both the pressure and convection terms as dominated by the
viscosity term. 
To begin to see why, recall that the Laplace operator
$\Delta\colon H^2(\Omg)\cap H^1_0(\Omg)\to L^2(\Omg)$ is an
isomorphism, and note that $\grad\ps$ is determined by $\lap\uu$ via
\begin{equation}\label{ps-proj}
\nabla \ps = (I-\PP - \QQ)\lap\uu, \qquad \QQ:=\nabla \div \Delta^{-1}.
\end{equation}
Equation \qref{NSE5} can then be written
\begin{eqnarray}
 \pa_t \uu + \PP( \udotgrad \uu - \ff)
 &=& \nu (\PP+\QQ) \lap\uu  \nonumber\\
 &=& \nu \lap \uu  - \nu(I-\PP-\QQ)\lap\uu.  \label{NSE7}
\end{eqnarray}
Theorem~\ref{T.main} will allow us to regard the last term as a controlled
perturbation.

This approach should be contrasted with the usual one that regards the
Navier-Stokes equations as a perturbation of the Stokes system.  
To show the advantage of this point of view, we will sketch the proof
from \cite{LLP} of {\em unconditional stability} of a simple
time-discretization scheme with explicit time-stepping for the
pressure and nonlinear convection terms and that is implicit only in
the viscosity term.

The discretization that we use is related to a class of extremely
efficient numerical methods for incompressible flow 
 \cite{Ti96,Pe,JL,GuS}.
Thanks to the explicit treatment of the convection and pressure
terms, the computation of the momentum equation is completely
decoupled from the computation of the kinematic pressure Poisson
equation used to enforce incompressibility.  No stationary Stokes
solver is necessary to handle implicitly differenced pressure
terms. For three-dimensional flow in a general domain, the
computation of incompressible Navier-Stokes dynamics is basically
reduced to solving a heat equation and a Poisson equation at each
time step. This class of methods is very flexible and can be used
with all kinds of spatial discretization methods \cite{JL},
including finite difference, spectral, and finite element methods.
The stability properties we establish here should be helpful in
analyzing these methods.

Indeed, we will show here that our stability analysis easily adapts to
proving unconditional stability for
corresponding fully discrete finite-element methods with $C^1$
elements for velocity and $C^0$ elements for pressure. 
(For additional details and convergence results see \cite{LLP}.) It is
important to note that we impose {\em no inf-sup compatibility
condition} between the finite-element spaces for velocity and
pressure.

Below, we will also describe how the unconstrained formulation above
can be extended to handle nonhomogeneous boundary conditions.  In
\cite{LLP} this formulation was exploited to establish a new result
regarding existence and uniqueness for strong solutions with minimal
compatibility conditions for total flux through the boundary.
For further results and more complete references and proofs, we refer to
\cite{LLP}.

\section{Properties of the Stokes pressure} \label{S.lessthan1}

Let $\Omg\subset\R^N$ be a bounded domain with $C^3$ boundary
$\Gam$. For any $\xx\in\Omg$ we let $\Phi(\xx) = \dist(x,\Gamma)$
denote the distance from $x$ to $\Gam$. For any $s>0$ we denote
the set of points in $\Omg$ within distance $s$ from $\Gamma$ by
\begin{equation} \label{E.Omgs}
\Omg_s = \{ \xx \in \Omg \mid \Phi(\xx) \leq s \},
\end{equation}
and set $\Omg_s^c = \Omg \backslash \Omg_s$ and $\Gam_s = \{ \xx
\in \Omg \mid \Phi(\xx) = s \}$. Since $\Gamma$ is $C^3$ and
compact, there exists $s_0>0$ such that $\Phi$ is $C^3$ in
$\Omg_{s_0}$ and its gradient is a unit vector, with $|\grad
\Phi(\xx)| = 1$ for every $\xx\in\Omg_{s_0}$. We let
\begin{equation} \label{E.nn}
   \nn(\xx) = -\grad \Phi(\xx),
\end{equation}
then $\nn(\xx)$ is the outward unit normal to $\Gam_s=\pa\Omg_s^c$
for $s=\Phi(\xx)$, and $\nn \in C^2(\bOmg_{s_0},\R^N)$.

We let $\<f,g\>_\Omg = \int_{\Omg}fg$ denote the $L^2$ inner product
of functions $f$ and $g$ in $\Omg$, and let $\| \cdot \|_\Omg$ denote
the corresponding norm in $L^2(\Omg)$.  We drop the subscript on the
inner product and norm when the domain of integration is understood in
context.

\subsection{Estimates related to the Neumann-to-Dirichlet map}
Taking the divergence and normal component of both side of \qref{StokesP}, 
we see that Stokes
pressure is harmonic in $\Omg$, and is determined as the zero-mean
solution of the Neumann boundary-value problem 
\begin{equation}\label{NeumannP}
\Delta\ps=0 \quad\mbox{in $\Omg$},\qquad
\nn\cdot\nabla\ps=\nn\cdot(\lap-\grad\div)\uu \quad\mbox{on $\Gamma$}.
\end{equation}
A useful ingredient in our proof of \qref{StokesPE} is an integrated 
version of a standard estimate that controls tangential gradients
at the boundary (these are determined from Dirichlet data) in
terms of the normal gradient (Neumann data). 
\begin{lemma} \label{L.N2D}
  Let $\Omg \subset \R^N$ be a bounded domain with $C^2$ boundary.
  Then there exist positive constants $r_0$, $C_1$ and
  $C_0$ such that for any $p$ that satisfies
  $\lap p =0$ in $\Omg$, for any $s \in(0, r_0]$ we have
   \begin{equation} \label{gradnder}
     \int_{\Gam_s} |\grad p |^2 \leq C_1 \int_{\Gam_s}
     |\nn \cdot \grad p|^2 ,
   \end{equation}
and furthermore,
   \begin{equation} \label{E.wkN2D}
     \int_\Omgs |\grad p|^2 \leq C_1
     \int_{\Omg_s} |\nder p |^2 .
   \end{equation}
\end{lemma}

\noind{\bf Proof:} 
If $\Omg$ has smooth boundary and  $s\in(0,s_0]$ is fixed, with $s_0$
taken as above, the estimate \qref{gradnder} is a consequence
of classical elliptic theory as developed in the book of Lions and
Magenes~\cite{LM}, applied to the Neumann problem
\begin{equation}\label{Neum}
\lap p = 0 \quad\text{ in }\Omg_s^c,
\qquad \nder p =g \quad\text{ on }\Gam_s.
\end{equation}
Under the solvability condition $\int_{\Gam_s}g=0$, this theory
yields a map $g \mapsto p \mapsto p|_{\Gam_s}$ from
$L^2(\Gam_s)\to H^{3/2}(\Omg_s)\to H^1(\Gam_s)$ that is bounded.
This bounds the tangential part of $\grad p$ in $L^2(\Gam_s)$, which
is sufficient to establish \qref{gradnder}.
We need to verify that these bounds are uniform in $s$ for domains
with $C^2$ boundary.  We defer details to an appendix, where we
sketch a streamlined version of the arguments of Lions and Magenes
\cite{LM} that yield the desired bounds, based on using trace theorems
and the Lax-Milgram
lemma to solve \qref{Neum} with $g\in H^{-1/2}(\Gamma_s)$, and
standard elliptic regularity theory to treat $g\in H^{1/2}(\Gamma_s)$,
then interpolating to handle $g\in L^2(\Gamma_s)$.
The key to obtaining uniform estimates is to work with a fixed system
of local boundary-flattening maps that work simultaneously for all
$\Gam_s$ --- these are constructed using the distance function $\Phi(x)$.

Once \qref{gradnder} is established, one obtains \qref{E.wkN2D}
by simply integrating \qref{gradnder} in $s$.

\subsection{Identities at the boundary}

A key part of the proof of Theorem~\ref{T.main} involves boundary
values of two quantities that involve the decomposition of
$\uu=(\tanproj)\uu+\nproj\uu$ into parts parallel and normal to
the boundary, for which we have the following Lemma. 
The proof involves only straightforward computations and a density
argument; please see \cite{LLP} for the details.
\begin{lemma} \label{L.uH2}
Let $\Omg \subset \R^N$ be a bounded domain with boundary
$\Gam$ of class $C^3$. Then for any $\uu \in H^2(\Omg, \R^N)$ with
$\uu|_\Gam =0 $, the following is valid on $\Gam$:
 \begin{itemize}
   \item[$(i)$] $\div \( (\tanproj)\uu \)=0 $ in $H^{1/2}(\Gam)$.
   \item[$(ii)$] $ \nn \cdot  ( \lap - \grad \div) \( \nproj \uu \) = 0 $
   in $H^{-1/2}(\Gam)$.
 \end{itemize}
\end{lemma}

\subsection{Identities for the Stokes pressure} \label{S.upeupa}

Given $\uu \in H^2(\Omg,\R^N) \cap H^1_0(\Omg,\R^N)$, using the
fact that $\PP \grad (\div \uu)=0$, from \qref{StokesP} we get
\begin{equation}\label{pdef}
\grad \ps = (I-\P) \lap \uu - \grad \div \uu =
(I-\P)(\lap-\grad\div)\uu
\end{equation}
Note that from the trace theorem for $\Hdiv$ \cite[Theorem
2.5]{GR}, whenever $\ba\in L^2(\Omg,\R^N)$ and $\div\ba=0$ and
$\nn\cdot\ba|_\Gam=0$, then $(I-\P)\ba=0$. Thus, the Stokes
pressure is not affected by any part of the velocity field that
contributes nothing to $\nn\cdot\ba|_\Gam$ where
$\ba=(\lap-\grad\div)\uu$. Indeed, this means that the Stokes
pressure is not affected by the part of the velocity field in the
interior of $\Omg$ away from the boundary, and it is not affected
by the normal component of velocity near the boundary, since
$\nn\cdot(\lap-\grad\div)(\nproj\uu)|_\Gam=0$ by
Lemma~\ref{L.uH2}.

This motivates us to focus on the part of velocity near and parallel
to the boundary. We make the following decomposition.
Fix $s>0$ so $2s<\min(s_0,r_0)$ and fix any smooth
non-negative cutoff function $\cut:\Omg\to[0,1]$ 
with $\cut=1$ on $\Omg_{s/2}=\{x\in\Omg\mid\dist(x,\Gam)\le \frac12 s\}$
and $\cut=0$ on $\Omgin$. Then we can write
\begin{equation}\label{upaupe}
\uu = \upe + \upa
\end{equation}
where
\begin{equation} \label{updefs}
     \upe = \cut \nproj \uu + (1-\cut) \uu , \qquad
     \upa = \cut (\tanproj)\uu  .
\end{equation}
Now $\upe = (\nproj) \uu$ in $\Omgbdry$, and with
$\ba_{\perp}=(\lap-\grad\div)\upe$ we have
\begin{equation}
\mbox{$\ba_{\perp}\in L^2(\Omg,\R^N)$,\quad
$\div\ba_{\perp}=0$,\quad and\ \ $\nn\cdot\ba_{\perp}|_\Gam=0$}
\end{equation}
by Lemma~\ref{L.uH2}(ii).  Hence $\<\ba_\perp,\grad\phi\>=0$ for
all $\phi\in H^1(\Omg)$, i.e.,
\begin{equation}\label{upezero}
(I-\P) (\lap - \grad \div )\upe=0.
\end{equation}
Combining this with \qref{upaupe} and \qref{pdef}
proves part (i) of the following.
\begin{lemma} \label{L.pair}
 Let $\Omg\subset\R^N$ be a bounded domain with $C^3$ boundary, and
 let $\uu \in H^2(\Omg,\R^N) \cap H_0^1(\Omg,\R^N)$.
Let $\ps$ and $\upa$ be defined as in
 \qref{pdef} and \qref{updefs} respectively. Then
 \begin{itemize}
   \item[$(i)$] The Stokes pressure is determined by $\upa$ according to
   the formula
 \begin{equation}\label{psupa}
 \grad \ps = (I-\P) (\lap - \grad \div) \upa .
 \end{equation}
   \item[$(ii)$] For any $q \in H^1(\Omg)$ that satisfies
   $\lap q =0$ in the sense of distributions,
             \begin{equation} \label{stokespair-1}
               \< \lap \upa -\grad \ps, \grad q \> = 0.
             \end{equation}
   \item[$(iii)$] In particular we can let $q=\ps$ in (ii), so
   $ \< \lap \upa-\grad\ps, \grad \ps \>  = 0$ and
\begin{equation} \label{stokespair-2}
\| \lap \upa \|^2 = \| \lap \upa - \grad \ps \|^2 + \| \grad \ps \|^2.
              \end{equation}
 \end{itemize}
\end{lemma}
{\bf Proof:} We already proved (i). For (ii), note
by Lemma~\ref{L.uH2}(i) we have
\begin{equation}
     \div \upa |_{\Gam} =0,
\end{equation}
so $\div\upa\in H^1_0(\Omg)$.
Using part (i), whenever $q\in H^1(\Omg)$ and $\lap q=0$ we get
\begin{equation}
\<\grad\ps,\grad q\> =
\<\lap\upa-\grad\div\upa, \grad q\> =
\<\lap\upa,\grad q\>.
\end{equation}
This proves (ii), and
then (iii) follows by the $L^2$ orthogonality. $\square$
\subsection{Proof of \qref{StokesPE} for {some} $\beta <1$}

Let $\uu\in H^2(\Omg,\R^N)\cap H^1_0(\Omg,\R^N)$ and define the
Stokes pressure $\grad\ps$ by \qref{StokesP} and the decomposition
$\uu=\upe+\upa$ as in the previous subsection. Then by part (iii) of
Lemma~\ref{L.pair} we have
\begin{equation}\label{lapu}
\|\lap\uu\|^2 = \|\lap\upe\|^2+ 2\<\lap\upe,\lap\upa\> +
\|\lap\upa-\grad\ps\|^2+\|\grad\ps\|^2.
\end{equation}
We will establish the Theorem with the help of two further estimates.

\noindent{\bf Claim 1:} For any $\eps>0$, there exists a constant
$C_\eps>0$ independent of $\uu$ such that
\begin{equation}\label{E.claim2}
\<\lap\upe,\lap\upa\> \ge -\eps\|\lap\uu\|^2 - C_\eps\|\grad\uu\|^2.
\end{equation}

\noindent{\bf Claim 2:} There exist positive constants $\pcon_1$ and
$C_2$ independent of $\uu$ such that
\begin{equation}
\|\lap\upa-\grad\ps\|^2 \ge \pcon_1 \|\grad\ps\|^2 -
C_2\|\grad\uu\|^2.
\end{equation}

\noind {\bf Proof of claim 1:} From the definitions in \qref{updefs},
we have
\begin{equation}\label{upests}
\lap\upe = \cut \nproj \lap \uu + (1-\cut)\lap \uu + R_1,
\qquad
\lap\upa =
    \cut (\tanproj) \lap \uu  + R_2,
\end{equation}
where $\|R_1\|+\|R_2\|\le C\|\grad\uu\|$ with $C$ independent of
$\uu$.
Since $\tanproj=(\tanproj)^2$,
\[
  \( \cut \nproj \lap \uu +(1-\cut)\lap \uu  \)\cdot
  \( \cut (\tanproj) \lap \uu \) =
  0+\cut (1-\cut) | (\tanproj) \lap \uu |^2 \geq 0.
\]
This means the leading term of $  \<\lap \upe ,\lap \upa\> $ is
non-negative. Using the inequality $|\<a,b\>|\le
(\eps/C)\|a\|^2+(4C/\eps)\|b\|^2$ and the bounds on $R_1$ and
$R_2$ to estimate the remaining terms, it is easy to obtain
\qref{E.claim2}.

\noindent{\bf Proof of claim 2:}
Recall that $\upa$ is supported in $\Omgs$, with
\begin{equation}
\lap\upa = \cut(\tanproj)\lap\uu + R_3
\end{equation}
where $\|R_3\|\le C\|\grad\uu\|$.
Since $\nn\cdot(\tanproj)\lap\uu=0$ we find 
\begin{equation}
\|\nn\cdot\lap\upa\|_{\Omgs} \le C_2 \|\grad\uu\|
\label{upaest00}
\end{equation}
with $C_2>0$ independent of $\uu$. 
Using $\|a+b\|^2+\|a\|^2 \ge
\frac12\|b\|^2$ and Lemma~\ref{L.N2D}, we get
\begin{align}
 \|\lap\upa-\grad\ps\|^2 \ge &\  %
 \int_{\Omg_s^c} |\grad\ps|^2 + 
 \int_\Omgs |\nn\cdot(\lap\upa-\grad\ps)|^2
 \nonumber\\ \ge &\ %
 \int_{\Omg_s^c} |\grad\ps|^2 + 
 \frac12 \int_\Omgs |\nn\cdot\grad\ps|^2 
 - \int_\Omgs |\nn\cdot\lap\upa|^2 
 \nonumber\\ \ge &\ %
 \pcon_1\|\grad\ps\|^2 - C_2 \|\grad\uu\|^2.
\label{upest0}
\end{align}
with $\pcon_1=\min(1,1/(2C_1))$. This establishes Claim 2.

Now we conclude the proof of the theorem. Combining the claims
with \qref{lapu}, we get
\begin{equation}
(1+2\eps) \|\lap\uu\|^2 \ge (1+\pcon_1)\|\grad\ps\|^2 -
(C_2+2C_\eps)\|\grad\uu\|^2.
\end{equation}
Taking $\eps>0$ so that $2\eps<\pcon_1$ yields \qref{StokesPE}
with $\pcon=(1+2\eps)/(1+\beta_1)<1$.
$\square$

\section{Unconditional stability of time discretization with pressure explicit}
\label{S.stab}

In this section we exploit Theorem~\ref{T.main} to establish the
unconditional stability of a simple time discretization scheme for the
initial-boundary-value problem for \qref{NSE5}, our unconstrained
formulation of the Navier-Stokes equations. We focus here on the case of
two and three dimensions.

Let $\Omg$ be a bounded domain in $\R^N$ ($N=2$ or 3) 
with boundary $\Gam$ of class $C^3$.
We consider the initial-boundary-value problem
 \begin{align}
 \pa_t \uu + \udotgrad \uu + \nabla \pe + \nu \nabla \ps
 = \nu \Delta \uu + \ff
  &\qquad (t>0,\ x\in\Omg), \label{newNSE1} \\
  \uu = 0
  &\qquad (t\ge0,\ x\in\Gam), \label{newNSE2}\\
  \uu=\uuin
  &\qquad (t=0,\ x\in\Omg). \label{newNSE3}
 \end{align}
We assume $\uuin \in H_0^1(\Omg, \R^N)$ and
$\ff \in L^2(0,T;L^2(\Omg, \R^N))$ for some given $T>0$.

Theorem~\ref{T.main} tells us that the Stokes pressure can be
strictly controlled by the viscosity term. This allows us to treat
the pressure term explicitly, so that the update of pressure is
decoupled from that of velocity. This can make corresponding
fully discrete numerical schemes very efficient (see \cite{JL},
also \cite{Ti96,Pe,GuS}).
Here, through Theorem~\ref{T.main}, we will prove that the following
spatially continuous time discretization scheme
has surprisingly good stability properties:
\begin{align}
  \frac{\uu^{n+1} - \uu^n }{\Dt} -  \nu \lap \uu^{n+1}
= \ff^n -\uu^n \cdot \grad \uu^n - \grad \pe^n - \nu \grad \ps^n ,
\label{fdiff1}  \\
    \grad \pe^n =  (I-\P) (\ff^n- \uu^n \cdot \grad \uu^n ) ,
\label{fdiff2-1}\\
    \grad \ps^n = (I-\P)\lap \uu^n - \grad (\div\uu^n) ,
\label{fdiff2-2}\\
  \uu^n\big|_{\Gam} = 0.      \label{fdiff3}
\end{align}
We set
\begin{equation}\label{fndef}
\ff^n = \frac{1}{\Dt} \int_{n\Dt}^{(n+1)\Dt} \ff(t)\, dt,
\end{equation}
and take $\uu^0\in \uspace$
to approximate $\uuin$ in $H^1_0(\Omg,\R^N)$.
It is evident that for all $n=0,1,2,\ldots$, given $\uu^n\in H^2\cap H^1_0$
one can determine $\grad\pe^n\in L^2$ and $\grad\ps^n\in L^2$
from \qref{fdiff2-1} and \qref{fdiff2-2}
and advance to time step $n+1$ by solving \qref{fdiff1} as an elliptic
boundary-value problem with Dirchlet boundary values to obtain
$\uu^{n+1}$.

Let us begin making estimates --- our main result is stated as
Theorem~\ref{T.stab} below.
Dot \qref{fdiff1} with $-\lap u^{n+1}$ and use \qref{fdiff2-1}
and $\|I-\PP\|\le1$ to obtain
\begin{align}
  \frac{1}{2\Dt} \Big( &\|\grad \uu^{n+1}\|^2  - \| \grad \uu^n\|^2 +
  \| \grad \uu^{n+1} - \grad \uu^n\|^2 \Big) + \nu \|\lap \uu^{n+1}\|^2
 \nonumber\\
  & \leq \|\lap \uu^{n+1} \| \Big(
  2\| \ff^n - \uu^n \cdot \grad \uu^n \|
 + \nu \|\grad \ps^n \| \Big)
\nonumber\\&
\leq \frac{\eps_1}{2}\|\lap\uu^{n+1}\|^2 + \frac{2}{\eps_1}
\|\ff^n-\uu^n\cdot\nabla\uu^n\|^2
+\frac{\nu}{2} \( \|\lap\uu^{n+1}\|^2+\|\grad\ps^n\|^2\)
\label{E.dot1}
\end{align}
for any $\eps_1>0$. 
This gives
\begin{align}
  \frac{1}{\Dt} \Big( & \| \grad \uu^{n+1}\|^2 - \| \grad \uu^n\|^2
  \Big) +
  (\nu-\eps_1) \|\lap \uu^{n+1}\|^2 \nonumber
 \\ & \le
\frac{8}{\eps_1}\left(
\|\ff^n\|^2+\|\uu^n\cdot\nabla\uu^n\|^2\right) +\nu \|\grad \ps^n
\|^2 .
  \label{E.temp.9}
\end{align}
By Theorem~\ref{T.main}, for some $\pcon<1$ one has
\begin{equation} \label{E.cite.them1}
 \nu \|\grad \ps^{n}\|^2 \le \nu\pcon \|\lap \uu^{n} \|^2
  +\nu C_{\pcon} \|\grad \uu^{n} \|^2.
\end{equation}
Using this in \qref{E.temp.9},  one obtains
\begin{align}
  \frac{1}{\Dt} \Big(  \|\grad \uu^{n+1}\|^2  -& \|\grad \uu^n\|^2
  \Big)
   + (\nu-\eps_1) \( \|\lap \uu^{n+1}\|^2 - \|\lap \uu^n\|^2 \)
   \nonumber\\
   &   + ( \nu - \eps_1 -\nu\pcon ) \|\lap \uu^n \|^2
    \nonumber \\
    &  \le \frac{8}{\eps_1}
 \( \|\ff^n\|^2 + \| \uu^n \cdot \grad \uu^n \|^2 \)
    + \nu C_{\pcon} \|\nabla \uu^n \|^2 .
    \label{E.fd-mid}
\end{align}

At this point the pressure has been dealt with.
Then, it is rather straightforward and
standard to derive from Ladyzhenskaya's inequalities (see \cite{LLP}) that
\begin{align} \label{E.nlin}
 \| \uu^n \cdot \grad \uu^n \|^2 & \leq
 \eps_2  \|\lap \uu^n \|^2 + C \| \grad \uu^n \|^6,
\end{align}
for any $\eps_2>0$. Plug this into \qref{E.fd-mid} and take $\eps_1$,
$\eps_2 >0$ satisfying $\nu - \eps_1 >0$ and
$\eps:=\nu - \eps_1-\nu\pcon  - 8 \eps_2 / \eps_1 > 0 $.
We get
\begin{align}
 \frac{1}{\Dt} \( & \|\grad \uu^{n+1}\|^2  - \|\grad \uu^n\|^2 \)
  + (\nu -\eps_1 ) \( \|\lap \uu^{n+1}\|^2 - \| \lap \uu^n \|^2 \)
  + \eps \| \lap \uu^n \|^2
    \nonumber \\
 &\leq   \frac{8}{\eps_1} \|\ff^n\|^2
   +
     C \|\grad \uu^{n} \|^6
   + \nu C_{\pcon} \|\grad \uu^{n}\|^2 . \label{fdest1}
\end{align}
A Gronwall-type argument now leads to 
a simplification of the stability result in \cite{LLP}:
\begin{theorem}\label{T.stab}
Let $\Omg$ be a bounded domain in $\R^N$ ($N=2$ or $3$) with $C^3$
boundary, and assume $\ff \in L^2(0,T;L^2(\Omg, \R^N))$ for some
given $T>0$ and $\uu^0 \in H_0^1(\Omg, \R^N)\cap H^2(\Omg,\R^N)$.
 Consider the time-discrete scheme \qref{fdiff1}-\qref{fndef}.
 Then there exist positive constants
 $T^*$ and $C_3$,
 such that whenever $n\Dt\le T^*$, we have
 \begin{align}
  \sup_{0\leq k \leq n }\|\grad \uu^k \|^2
  + \sum_{k=0}^n  \|\lap \uu^k \|^2 \Dt \le C_3 .
  \label{stab1} 
 \end{align}
The constants $T^*$ and $C_3$ depend only upon $\Omg$, $\nu$ and
\[
M_0:=\|\grad\uu^0\|^2+\nu\Dt\|\lap\uu^0\|^2+\int_0^T\|\ff\|^2.
\]
\end{theorem}

\noindent{\bf Proof:} Put
\begin{equation}
z_n = \| \grad \uu^n \|^2 + (\nu-\eps_1) \Dt \|\lap \uu^n\|^2 ,
\quad w_n = 
\eps \|\lap\uu^{n}\|^2 , \quad b_n = \|\ff^n\|^2,
\end{equation}
and note that from \qref{fndef} we have that 
$\sum_{k=0}^{n-1}b_k^2\Dt \le \int_0^T |\ff(t)|^2\,dt $
as long as $n\Dt\le T$.
Then by \qref{fdest1},
\begin{equation}
z_{n+1} +w_n\Dt \le z_n+ C\Dt(b_n+ z_n + z_n^3) , \label{zeq1}
\end{equation}
where we have replaced $\max \{ 8/\eps_1, C, 
\nu C_{\pcon}  \} $ by $C$. Summing from 0 to $n-1$ 
yields
\begin{equation}
z_n+\sum_{k=0}^{n-1}w_k\Dt \le CM_0 + C\Dt\sum_{k=0}^{n-1}(z_k+z_k^3)
=: y_n.
\label{ydef}
\end{equation}
The quantities $y_n$ so defined increase with $n$ and satisfy
\begin{equation}
y_{n+1}-y_n = C\Dt(z_n+z_n^3) \le C\Dt(y_n+y_n^3).
\label{yeq1}
\end{equation}
Now set $F(y)=\ln(\sqrt{1+y^2}/y)$ so that
$F'(y)=-(y+y^3)^{-1}$. Then on $(0,\infty)$,
$F$ is positive, decreasing and convex, and we have
\begin{equation}
F(y_{n+1})-F(y_n) = F'(\xi_n)(y_{n+1}-y_n) \ge
- \frac{y_{n+1}-y_n}{y_n+y_n^3} \ge -C\Dt,
\end{equation}
whence
\begin{equation}
F(y_n)\ge F(y_0)- Cn\Dt= F(CM_0)-Cn\Dt.
\end{equation}
Choosing any $T^*>0$ so that $C_*:= F(CM_0)-C T^*>0$, we infer that
as long as $n\Dt\le T^*$ we have $y_n\le F^{-1}(C_*)$, and this
together with \qref{ydef}
yields the stability estimate \qref{stab1}.
$\square$

\section{Unconditional stability for $C^1/C^0$
finite element methods without inf-sup conditions}
\label{S.fem}

The simplicity of the stability proof for the time-discrete scheme
above allows us to easily establish the unconditional stability 
of corresponding fully discrete finite-element methods that use $C^1$
elements for the velocity field and $C^0$ elements for pressure. 

We suppose that for some sequence of positive values of $h$
approaching zero, $\xuh\subset\uspace$ is a finite-dimensional space
containing the approximate velocity field, and suppose $\xph\subset
H^1(\Omg)/\R$ is a finite-dimensional space containing approximate
pressures. 
We discretize \qref{newNSE1} in a straightforward
way, implicitly only in the viscosity term and explicitly in the
pressure and nonlinear terms. The resulting scheme was also derived 
in \cite{JL} and is equivalent to a
space discretization of the scheme in \qref{fdiff1}--\qref{fndef}.
Given the approximate velocity
$\uu^h_n$ at the $n$-th time step, we determine 
$p_h^n\in\xph$ and $\uu_h^{n+1}\in\xuh$ by requiring
\begin{align} 
\<\grad p_h^n +\nu\grad\div\uu_h^n-\nu\lap\uu_h^n
+\uu_h^n\!\cdot\!\grad\uu_h^n-\ff^n,\grad \phi_h\>=0
\quad &\forall \phi_h \in \xph, \label{fem.1}
  \\
  \< \frac{\grad\uu^{n+1}_h - \grad\uu^n_h }{\Dt}, \grad \vv_h \>  + 
\<\nu\lap\uu_h^{n+1},\lap\vv_h\>
= \<\grad p_h^n 
+\uu_h^n\!\cdot\!\grad\uu_h^n & -\ff^n,\lap\vv_h \>
\nonumber \\ 
  & \forall \vv_h \in \xuh . 
\label{fem.2} 
\end{align}

We are to show the scheme above is unconditionally stable. 
First, we take $\phi_h = p_h$ in \qref{fem.1}. Due to the fact that
$ \<\PP(\lap-\grad\div)\uu_h^n,\grad p_h^n\>=0, $
we directly deduce from the Cauchy-Schwarz inequality that
\begin{equation}\label{E.phest}
  \|\grad p_h^n \| \le
  \|\nu\grad\ps(u_h^n)\|+\|\uu_h^n\!\cdot\!\grad\uu_h^n - \ff^n\|
\end{equation}
where
\begin{equation}  
  \grad\ps(u_h^n) = (I-\PP)(\lap-\grad\div)\uu_h^n
\end{equation}
is the Stokes pressure associated with $\uu_h^n$. (Note
$\grad\ps(u^n_h)$ need not lie in the space $\xph$). 
Now, taking $\vv_h = \uu_h^{n+1}$ in \qref{fem.2} and arguing just
as in \qref{E.dot1}, 
we obtain an exact analog of \qref{E.temp.9}, namely
\begin{align}
  \frac{1}{\Dt} \Big( & \| \grad \uu_h^{n+1}\|^2 - \| \grad \uu_h^n\|^2
  \Big) +
  (\nu-\eps_1) \|\lap \uu_h^{n+1}\|^2 \nonumber
 \\ & \le
\frac{8}{\eps_1}\left( \|\ff^n\|^2+\|\uu_h^n\cdot\nabla\uu_h^n\|^2\right)
+\nu \|\grad \ps(\uu_h^n) \|^2 .
  \label{fem.ineq}
\end{align}
Proceeding now exactly as in section 3 leads to the following result. 

\begin{theorem}\label{T.stabh}
Let $\Omg$ be a bounded domain in $\R^N$ ($N=2$ or $3$) with $C^3$
boundary, and suppose $\xuh\subset\uspace$, $\xph\subset
H^1(\Omg)/\R$.  Assume $\ff \in L^2(0,T;L^2(\Omg, \R^N))$ for some
given $T>0$ and $\uu_h^0 \in \xuh$.
 Consider the finite-element scheme \qref{fem.1}-\qref{fem.2} with
 \qref{fndef}.  Then there exist positive constants
 $T^*$ and $C_4$,
 such that whenever $n\Dt\le T^*$, we have
 \begin{align}
  \sup_{0\leq k \leq n }\|\grad \uu^k_h \|^2
  + \sum_{k=0}^n  \|\lap \uu^k_h \|^2 \Dt \le C_4 .
  \label{stab1h} 
\end{align}
The constants $T^*$ and $C_4$ depend only upon $\Omg$, $\nu$ and 
\[
M_{0h}:=\|\grad\uu_h^0\|^2+\nu\Dt\|\lap\uu_h^0\|^2+\int_0^T\|\ff\|^2.
\]
\end{theorem}

In \cite{LLP} further arguments are given that establish the 
convergence of these approximation schemes 
up to the maximal time of existence for the strong solution of
the unconstrained Navier-Stokes formulation
\qref{newNSE1}-\qref{newNSE3}.

\section{Non-homogeneous side conditions} \label{S.nonhom}
Consider the Navier-Stokes equations with non-homogeneous boundary
conditions and divergence constraint:
 \begin{align}
  \pa_t \uu + \uu \cdot\! \grad \uu + \nabla p  =  \nu \Delta \uu
  + \ff  & \qquad (t > 0, x \in \Omg ),\label{nh-oldNSE1} \\
  \nabla \cdot \uu  =  h  & \qquad (t \geq 0, x \in \Omg ),
  \label{nh-oldNSE2} \\
  \uu   =  \gg  & \qquad (t \geq 0, x \in \Gam ), \label{nh-oldNSE3} \\
  \uu  = \uuin  & \qquad (t = 0, x \in \Omg ). \label{nh-oldNSE4}
\end{align}
What we have done before can be viewed as replacing the
divergence constraint \qref{nh-oldNSE2} by decomposing the pressure
via the formulae in \qref{EulerP} and \qref{StokesP} in such a way
that the divergence constraint is enforced automatically.
It turns out that in the non-homogeneous case a very similar procedure
works. One can simply use the Helmholtz decomposition to identify
Euler and Stokes pressure terms {\it exactly as before} via the
formulae \qref{EulerP} and \qref{StokesP}, but in addition another
term is needed in the total pressure to deal with the inhomogeneities.
Equation \qref{NSE5} is replaced by
\begin{equation}\label{E.alt2}
 \pa_t \uu + \PP( \udotgrad \uu - \ff - \nu \Delta \uu ) +\grad\pgh
 = \nu \nabla(\div \uu).
\end{equation}
The equation that determines the inhomogeneous pressure $\pgh$ can be
found by dotting with $\grad \phi$ for $\phi\in H^1(\Omg)$,
formally integrating by parts and plugging in the side conditions:
We require
\begin{equation}\label{pgh1}
\<\grad\pgh,\grad \phi\> = -\<\pa_t (\nn \cdot \gg),\phi\>_\Gamma
+\<\pa_t h,\phi\> + \<\nu\grad h,\grad \phi\>
\end{equation}
for all $\phi\in H^1(\Omg)$. With this definition, we see from
\qref{E.alt2} that
\begin{equation}\label{wk-heat2}
\<\pa_t\uu,\grad \phi\>-\<\pa_t(\nn\cdot \gg),
\phi\>_\Gamma+\<\pa_th,\phi\> =\<\nu\grad(\div\uu-h),\grad \phi\>
\end{equation}
for every $\phi\in H^1(\Omg)$.
This will mean $w:=\div \uu -h $ is a weak solution of
\begin{equation} \label{wk-heat4}
  \pa_t w  = \nu \lap w \; \text{  in } \Omg, \qquad
  \nn \cdot\! \grad w  =  0 \; \text{  on } \Gam,
\end{equation}
with initial condition $w=\div \uuin -h\big|_{t=0}$. So the
divergence constraint will be enforced through exponential
diffusive decay as before (see \qref{nh-diss.iden} below).

The total pressure in \qref{nh-oldNSE1} now has the representation
 \begin{equation} \label{nh-ptotal}
   p = \pe + \nu \ps + \pgh,
 \end{equation}
where the Euler pressure $\pe$ and the Stokes pressure $\ps$ are
determined exactly by \qref{EulerP} and \qref{StokesP} as before, and
$\pgh$ is determined up to a constant by the forcing functions $g$ and
$h$ through the weak-form pressure Poisson equation \qref{pgh1}.
Our unconstrained formulation of
\qref{nh-oldNSE1}-\qref{nh-oldNSE4} then takes the form
 \begin{align}
 \pa_t \uu + \udotgrad \uu +  \nabla \pe +\nu \nabla \ps+ \grad
 \pgh = \nu \Delta \uu + \ff
  &\qquad (t>0,\ x\in\Omg), \label{nh-newNSE1} \\
  \uu = \gg
  &\qquad (t\ge0,\ x\in\Gam), \label{nh-newNSE2}\\
  \uu=\uuin
  &\qquad (t=0,\ x\in\Omg). \label{nh-newNSE3}
 \end{align}

We shall state an existence and uniqueness result for strong
solutions of the unconstrained formulation
\qref{nh-newNSE1}--\qref{nh-newNSE3}. We refer to \cite{LLP} for
the proof, which is based on using a classical trace theorem of Lions
and Magenes to reduce the problem to one with homogeneous boundary and
initial conditions. Then one uses the stability of the 
time-differencing scheme together with a standard weak compactness
method to get the existence.

Let $\Omg$ be a bounded, connected domain in $\R^N$ ($N=2$ or $3$)
with boundary $\Gam$ of class $C^3$. We assume
\begin{align}
 \uuin &\in 
                        H^1(\Omg,\R^N) ,  \label{nh-conduin} \\
 \ff &\in 
                  L^2(0,T;L^2(\Omg,\R^N)) , \label{nh-condf} \\
 \gg &\in 
                  H^{3/4}(0,T;L^2(\Gam,\R^N)) \cap L^2(0,T;
H^{3/2}(\Gam,\R^N)) \nonumber \\
& \hspace{1.7cm} \cap \{ \gg \; \big| \; \pa_t (\nn \cdot \gg)
\in L^2(0,T;H^{-1/2}(\Gam)) \}  , \label{nh-condg} \\
 h &\in H_h :=
                L^2(0,T;H^1(\Omg))   \cap H^1(0,T;(H^1)'(\Omg)) .
\label{nh-condh}
 \end{align}
Here $(H^1)'$ is the space dual to $H^1$. We also make the
compatibility assumptions
\begin{align}
&   \gg=\uuin \quad \text{ when $t=0$, $x \in \Gam$},
      \label{nh-compat} \\
&   \<  \pa_t (\nn \cdot \gg) , 1 \>_{\Gam}=\< \pa_t h, 1
\>_{\Omg}.
      \label{nh-conddivu}
\end{align}

We define
\begin{align}
& V := L^2(0,T;H^2(\Omg,\R^N))\cap H^1(0,T;L^2(\Omg,\R^N)) ,
\label{Vdef2} 
\end{align}
and note we have the embeddings (\cite[p.~288]{Ev},
\cite[p.~176]{Te})
\begin{equation}\label{VWemb}
V\emb C([0,T],H^1(\Omg,\R^N)),\quad H_h\emb C([0,T],L^2(\Omg)).
\end{equation}

\begin{theorem}\label{T.nonhomog}
Let $\Omg$ be a bounded, connected domain in $\R^N$ ($N=2$ or $3)$ and
assume \qref{nh-conduin}-\qref{nh-conddivu}. Then
there exists $T^*>0$ so that a unique strong solution of
\qref{nh-newNSE1}-\qref{nh-newNSE3} exists on $[0,T^*]$, with
\begin{align*}
\uu\in&
L^2(0,T^*;H^2(\Omg,\R^N))\cap H^1(0,T^*;L^2(\Omg,\R^N)),\\
p=&\nu\ps+ \pe + \pgh \in L^2(0,T^*;H^1(\Omg)/\R) ,
\end{align*}
where $\pe$ and $\ps$ are defined in \qref{EulerP} and
\qref{StokesP}, $\pgh \in L^2(H^1(\Omg)/\R)$ satisfies
\qref{pgh1}. 

Moreover, $\uu\in C([0,T^*],H^1(\Omg,\R^N))$ and
$$ \div \uu-h \in L^2(0,T^*;H^1(\Omg)) \cap
H^1(0,T^*;(H^1)'(\Omg))
$$
is a smooth solution of
the heat equation for $t>0$ with no-flux boundary conditions.
The map $t\mapsto\|\div\uu-h\|^2$ is smooth for $t>0$
and we have the dissipation identity
\begin{equation} \label{nh-diss.iden}
  \frac{d}{dt}\frac12\|\div \uu-h \|^2 +  \nu \|\grad (\div \uu -h)\|^2 =0.
\end{equation}
If we further assume
$h \in L^2(0,T;H^2(\Omg))\cap H^1(0,T;L^2(\Omg))$ and $\div \uuin \in H^1(\Omg)$,
then
\[
\div \uu \in L^2(0,T^*;H^2(\Omg)) \cap H^1(0,T^*;L^2(\Omg)).
\]
\end{theorem}
\bigskip

\section*{Acknowledgments}
The fact (related to Lemma~\ref{L.N2D})
that $\nder p\in L^2(\Omg_s)$ implies $\grad p\in L^2(\Omg)$
for harmonic $p$ was proved some years ago by Oscar Gonzalez and RLP
(unpublished) through a partitioning and flattening argument.
RLP is grateful for this collaboration.
This material is based upon work supported by the National Science
Foundation under grant no.\ DMS 03-05985 (RLP) and DMS-0107218 (JGL).


\appendix
\section{Uniform bounds on the Neumann-to-Dirichlet map}

Here we provide a more detailed sketch of the proof of inequality
\eqref{gradnder} in Lemma~\ref{L.N2D}, based on summarizing the
relevant arguments from \cite[vol.~I]{LM} and taking into account
the finite regularity of the boundary.

Without loss of generality we may assume $\Omg$ is connected,
so the same holds for $\Omg_s$ for small $s$ (see below).
For solutions of the Neumann problem
\begin{equation}\label{a.Neum}
\lap p = 0 \quad\text{ in }\Omg_s^c,
\qquad \nder p =g \quad\text{ on }\Gam_s,
\end{equation}
we seek to bound the tangential gradient of $p$ by $g$ in
$L^2(\Gam_s)$, uniformly for small $s$.
This is accomplished by interpolating between maps
\begin{align*}
 H^{-1/2}(\Gam_s) &\stackrel{A_1}{\to} H^1(\Omg_s)
 \stackrel{T_1}{\to} H^{1/2}(\Gam_s), \\
 H^{1/2}(\Gam_s) &\stackrel{A_2}{\to} H^2(\Omg_s)
 \stackrel{T_2}{\to} H^{3/2}(\Gam_s) .
\end{align*}
Here the maps $T_1$ and $T_2$ are the trace maps.
The map $A_1$ gives the zero-average weak solution of the Neumann
problem \qref{Neum}, by applying the Lax-Milgram lemma to the weak
form of \qref{a.Neum}:
\begin{equation}
\int_{\Omg_s^c}\grad p\cdot \grad q = \int_{\Gam_s}g q
\qquad\mbox{for all $q\in H^1(\Omg_s^c)$ with}\ \int_{\Omg_s^c}q=0.
\end{equation}
For our given solution of \qref{a.Neum}
we have $p-\pavs=A_1g$ where $\pavs = |\Omg_s^c|^{-1}\int_{\Omg_s^c}p$.
Since
\begin{equation}\label{a.est1}
\left| \int_{\Omg_s^c}gq \right| \le
\|g\|_{H^{-1/2}(\Gam_s)}\|T_1 q\|_{H^{1/2}(\Gam_s)} \le
C_0 \|g\|_{H^{-1/2}(\Gam_s)}\|q\|_{H^1(\Omg_s^c)}
\end{equation}
by a trace theorem \cite[vol.~I, p.~41]{LM}, the bound
\begin{equation}\label{a.pavbd}
\|p-\pavs\|_{H^1(\Omg_s^c)}\le C_1\|g\|_{H^{-1/2}(\Gam_s)}
\end{equation}
follows by taking $q=p-\pavs$ and using Poincar\'e's inequality.

The map $A_2$ is the restriction of $A_1$ and is bounded by
elliptic regularity theory
using only the $C^2$ regularity of the boundary.
Using that the trace maps $T_1$ and $T_2$ are bounded,
we obtain bounds
\begin{equation}\label{Tbounds}
\|p-\pavs\|_{H^{1/2}(\Gam_s)}\le C_2 \|g\|_{H^{-1/2}(\Gam_s)}, \qquad
\|p-\pavs\|_{H^{3/2}(\Gam_s)}\le C_3 \|g\|_{H^{1/2}(\Gam_s)}.
\end{equation}
By an abstract interpolation theorem \cite[vol.~I, p.~27]{LM}, the
map $T_1A_1$ restricts to a bounded map between interpolation
spaces,
\begin{equation}
TA: [H^{-1/2}(\Gam_s), H^{1/2}(\Gam_s)]_{1/2} \to
[H^{1/2}(\Gam_s), H^{3/2}(\Gam_s)]_{1/2}
\end{equation}
with $\|TA\|\le C_4$ where $C_4$ depends only on $C_2$ and $C_3$.
By \cite[vol.~I, p.~36]{LM} we have the isomorphisms
\begin{equation}
L^2(\Gam_s)\cong [H^{1/2}(\Gam_s), H^{3/2}(\Gam_s)]_{1/2}, \qquad
H^1(\Gam_s)\cong [H^{-1/2}(\Gam_s), H^{1/2}(\Gam_s)]_{1/2}
\end{equation}
with equivalent norms, and so the map $g\mapsto TAg = (p-\pavs)|_{\Gam_s}$
is bounded from $L^2(\Gam_s)$ to $H^1(\Gam_s)$.  This leads to the
estimate
\begin{align}
\|\grad p\|^2_{L^2(\Gam_s)} &= \|\grad (p-\pavs)\|^2_{L^2(\Gam_s)}
\nonumber\\
&  \leq \|p-\pavs\|^2_{H^1(\Gam_s)} + \|\nder p\|^2_{L^2(\Gam_s)}
  \leq C_4 \|\nder p \|^2_{L^2(\Gam_s)},
\end{align}
which corresponds to \qref{gradnder}.

We now argue that these bounds are uniform in $s$ for $0<s\le s_0$ small.
For this purpose one should consider how the fractional-order Sobolev
spaces $H^r(\Gam_s)$ and the norms in those spaces are defined, and
how one obtains the bounds in trace theorems, in Poincar\'e's
inequality, in the elliptic regularity estimate, and in comparing
equivalent norms.
All this can be done (except for Poincar\'e's inequality, which we
discuss below) using a fixed partition of unity on $\Omg$ and a fixed
family of maps that locally flatten the boundaries $\Gamma_s$
simultaneously.

To describe how this can work, let $\Phi(x)=\pm\dist(x,\Gam)$ be the
signed distance function on $\R^N$, positive in $\Omg$ and negative
outside.  Then for some $s_0>0$, $\Phi$ is $C^2$ for $0\le |\Phi(x)|\le s_0$.
Let $x_0\in\Gam$ be arbitrary. Translate and rotate coordinates as necessary so
$x_0=0$ and the normal $\nn(x_0)=-\vec{e}_N$.
For $x=(\hat{x},x_N)\in \R^{N-1}\times\R$, let
\begin{equation}\label{psidef}
\psi(x)=(\hat{x},\Phi(\hat{x},x_N)).
\end{equation}
Evidently the map $\psi$ locally flattens all $\Gamma_s$
simultaneously --- we have $\psi(x)=(\hat{x},s)$ for all
$x\in\Gam_s$.  By the proof of the inverse function
theorem, $\psi$ is a $C^2$ diffeomorphism of $B(x_0,r_0)$ onto its
image, for some $r_0\in(0,s_0]$ fixed and independent of $x_0$.

Let $r_k=2^{-k}r_0$, $k=1,2,\ldots$.
The compact set $\bar\Omg_{r_2}$ is covered by the union of balls
$B(x_0,r_1)$ for $x_0\in\Gam$. Hence there is a finite subcover
by balls $B(x^{(j)},r_1)$, $j=1,\ldots,m$. We denote by
$\psi_j$ the associated diffeomorphism on $B_j=B(x^{(j)},r_0)$; we
have $\vec{e}_N\cdot\psi_j(x)=\Phi(x)$ for $|x-x^{(j)}|<r_0$.
We fix a partition of unity on $\bar\Omg_{r_2}$ subordinate to this
finite cover by $B_j$. I.e., we fix $\alpha_j\in C^\infty_0(B_j)$,
$j=1,\ldots,m$, such that
$\sum_{j=1}^m \alpha_j(x)=1$ for all $x\in\bar\Omg_{r_2}$.

The norms in the fractional-order Sobolev spaces $H^r(\Gam_s)$
can then be defined as follows for $|r|\le2$.
Given $s\in[0,r_2]$, the norm of any $u\in C^2(\Gam_s)$ in
$H^r(\Gam_s)$ is given by  the norm of the $C^2$ function
$\hat x\mapsto \{ (\alpha_j u)\circ \psi_j^{-1}(\hat x,s)\}_{j=1}^m$
in $H^r(\R^{N-1})^m$. With this fixed system of cutoffs and flattening
maps, one can check that
the trace theorem for $T_1$ and $T_2$, the elliptic regularity
theory for $A_2$, and the norm equivalences that must be invoked to
obtain the estimates in \qref{Tbounds}, all involve bounds that are
valid uniformly for $s\in(0,r_1]$. It remains to check Poincar\'e's
inequality.

\subsection*{Uniform bounds in Poincar\'e's inequality.}
In the procedure that we have outlined, to obtain a uniform bound in
\qref{a.pavbd} we need that the constant in Poincar\'e's inequality
\begin{equation}\label{uni-p}
\int_{\Omg_s^c}
|p-\pavs|^2 \le
C \int_{\Omg_s^c}
|\grad p|^2, \qquad \pavs = \frac1{|\Omg_s^c|}\int_{\Omg_s^c}p,
\end{equation}
can be taken independent of $s$ for small $s$. Showing this is
somewhat related
to showing that $\Omg_s^c$ is connected for small $s$, which we can do
as follows: Since $\Omg$ is connected, any two points in $\Omg_s^c$
can be connected by a path in $\Omg$. On this path we can replace
points in $\Omg_s$ by the projection along normals onto $\Gam_{s'}$,
where $s'-s$ is positive and so small that the new path is contained in
$\Omg_s^c$ and still connects the same two points.

To prove the uniform Poincar\'e inequality \qref{uni-p}, we reduce to a
fixed $s$ by integrating along normals as follows.
Recall that the distance function
$\Phi(x)=\dist(x,\Gamma)$ is $C^2$ for $\Phi(x)\le s_0$.
Fix $\epss=\frac12 s_0$, and suppose $0<s\le \epss$.
Let $q$ be a smooth function on $\Omg_s^c$.
Then if $s<\Phi(x)<\epss$ we integrate
$\pa_t q(x-t\nn(x))$ from $t=0$ to $\epss$ and apply the
Cauchy-Schwarz inequality to get
\begin{equation}
|q(x)|^2 \le 2|q(x-\epss\nn(x))|^2 + 2\epss
\int_0^\epss |\grad q(x-t\nn(x))|^2 \,dt .
\end{equation}
Integrating this in the domain where $s<\Phi(x)<\epss$ and changing
variables on the right hand side we get
\begin{equation}
\int_{s<\Phi(x)<\epss} |q(x)|^2\,dx \le C
\left( \int_{\epss<\Phi(y)<2\epss}|q(y)|^2\,dy +
\int_{s<\Phi(y)<2\epss}|\grad q(y)|^2\,dy\right)
\end{equation}
with $C$ independent of $q$ and $s$.
Now, given a smooth function $p$ on $\Omg_s^c$, use this inequality
for $q=p-\bar{p}_\epss$, together with the standard Poincar\'e inequality
for the domain $\Omg_\epss^c$, and the $L^2$-optimality of the average
in $\Omg_s^c$, to conclude that
\begin{equation}
\int_{\Omg_s^c}|p-\pavs|^2 \le
\int_{\Omg_s^c}|p-\bar{p}_\epss|^2 \le
C\left( \int_{\Omg_\epss^c} |p-\bar{p}_\epss|^2 + \int_{\Omg_s^c}|\grad p|^2
\right)
\le C \int_{\Omg_s^c}|\grad p|^2
\end{equation}
for some $C$ independent of $p$ and $s$.



\end{document}